\documentclass[12pt]{article}
\usepackage{latexsym}
\def\rref#1{(\ref{#1})}
\def\be{\begin{equation}}
\def\ee{\end{equation}}
\def\bea{\begin{eqnarray}}
\def\eea{\end{eqnarray}}
\def\beas{\begin{eqnarray*}}
\def\eeas{\end{eqnarray*}}

\def \b {\beta}
\def \a {\alpha}
\def \g {\gamma}

\def\lR{{\hbox{{\rm I}\kern-.2em\hbox{\rm R}}}}
\def\IZ{{\hbox{{\rm I}\kern-.2em\hbox{\rm Z}}}}
\def\IC{{\hbox{{\rm I}\kern-.2em\hbox{\rm C}}}}
\def\II{{\hbox{{\rm I}\kern-.2em\hbox{\rm I}}}}

\begin{document}

\begin{titlepage}
\vspace{.5in}
\begin{flushright}
DFTT/60/99\\
IST/DM/31/99\\
November 1999\\
math.QA/9911015\\
\end{flushright}
\begin{center}
{\Large\bf
Quantum Matrix Pairs}\\
\vspace{.4in}
{J.\ E.~N{\sc elson}\footnote{\it email: nelson@to.infn.it}\\
       {\small\it Dipartimento di Fisica Teorica }\\
       {\small\it Universit\`a degli Studi di Torino}\\
       {\small\it via Pietro Giuria 1, 10125 Torino}\\
       {\small\it Italy}}\\
\vspace{1ex}
{\small and}\\
\vspace{1ex}
{R.\ F.~P{\sc icken}\footnote{\it  email:  picken@math.ist.utl.pt}\\
{\small\it Departamento de Matem\'{a}tica and Centro de
Matem\'{a}tica Aplicada}\\
{\small\it  Instituto Superior T\'{e}cnico}\\
{\small\it Avenida Rovisco Pais, 1049-001 Lisboa}\\
{\small\it Portugal}}
\\\vglue 5mm
\end{center}
\begin{center}
\begin{minipage}{4.8in}
\begin{center}
{\large\bf Abstract}
\end{center}
{\small
The notion of quantum matrix pairs is defined. These are pairs of 
matrices with non-commuting entries, which have  
the same pattern of 
internal relations,  
q-commute with each other under matrix multiplication, and  
are such that
products of powers of the matrices obey the same pattern of  internal 
relations as the original pair. 
Such matrices appear in an approach by the authors to quantizing gravity 
in 2 space 
and 1 time dimensions with negative cosmological constant
on the torus.
Explicit examples and 
transformations which generate new pairs from a given pair are presented. }
\end{minipage}
\end{center}
\end{titlepage}
\addtocounter{footnote}{-2}

\section{Introduction \label{sec1}}
The notion of quantum matrix pairs arose in the context of our recent work
\cite{2+1} on
quantum gravity in 2 space 
and 1 time dimensions with negative cosmological constant
on the torus. As shown by Witten \cite{wit}, this model is equivalent to a
Chern-Simons theory with 
non-compact structure group $SL(2, R)\times SL(2, R)$. After imposing the
constraints the classical 
geometry may be encoded, up
to equivalence, by
two  commuting $SL(2,R)$ matrices $U_1$ and $U_2$ (together with an identical
second pair for the 
other $SL(2,R)$ factor of the structure group), which represent the holonomies 
of a flat connection around two generating
cycles of the (abelian)
fundamental group of the torus. The usual approach to quantizing this theory,
e.g. \cite{NRZ}, 
is to work with gauge-invariant variables, namely the traces of the holonomies,
but in \cite{2+1} we 
chose instead to use the gauge covariant holonomy matrices themselves as
variables. There we also 
argued that, after quantization, these matrices obey a q-commutation relation
\be
U_1U_2=qU_2U_1,
\label{fund}
\ee
where $q\rightarrow 1$ in the classical limit. Using the gauge covariance to put
the matrices in 
standard form, the simplest solution of \rref{fund} is to take both matrices to
be diagonal, in which 
case the diagonal elements of $U_1$ and $U_2$ obey standard q-commutation
relations. However, 
when studying more general solutions with both matrices of upper-triangular form,
we found matrices 
which had non-trivial internal relations, in addition to the ``mutual'' relations
involving elements of 
both matrices. This feature of having internal relations is characteristic of
quantum groups, and indeed 
the algebraic structure which emerges has many similarities with quantum groups,
as we will explain 
below. The main purpose of this article is to describe this new algebraic
structure, and present some 
examples.
 
As mathematical objects quantum matrix pairs may be thought of
as a simultaneous  
generalization of two familiar notions of ``quantum mathematics'',  
namely the quantum  
plane and quantum groups. 
\footnote{For background material on noncommutative geometry and quantum groups,
see 
\cite{connes}, \cite{charipressley}, \cite{kassel}. } 
To make this  
statement more precise we first  
recall the quantum plane (over the field $k$), described by two 
non-commuting coordinates $x$ and $y$, satisfying the relation 
\be 
xy=qyx 
\label{qplane} 
\ee 
for some invertible $q\in k$, $q\neq 1$. The algebra of polynomial functions on
the 
quantum  
plane is then given by $k\{x,y\}/(xy-qyx)$, where $k\{x,y\}$ is the free 
algebra with coefficients in $k$ and  
$(xy-qyx)$ is the ideal generated by $xy-qyx$. This algebra is  a  
deformation of the algebra of  
polynomial functions in two commuting variables $k[x,y]$. 
 
Quantum groups may be presented in a similar way.
For example, consider a  
$2\times 2$ matrix  
\[ 
U=\left( \begin{array}{cc}a &b\\ c&d\end{array} \right) 
\label{qmatrix} 
\] 
 with non-commuting entries satisfying the  relations (R) 
\beas  
ab=qba; \quad &ac=qca;\quad &ad-da=(q-q^{-1})bc; \\  
      bc=cb; \quad &bd=qdb; \quad &cd=qdc. 
\label{qgrouprelns} 
\eeas 
for some invertible $q\in k$, $q\neq 1$. The algebra of polynomial functions of
these entries,  
denoted $M_q(2)$, is 
given by  
the quotient of the free algebra $k\{a,b,c,d\}$ by the ideal generated by 
the  
relations $(R)$, and is a deformation of the algebra of polynomial  
functions in the commuting variables $a,b,c,d$.   
 
The pattern of relations (R) seems  somewhat arbitrary at first sight, although
there are deep reasons for 
this particular form. One might also ask why the symbols $a,b,c,d$ are displayed
as entries of a matrix, 
instead of, say, as components of a $4$-vector belonging to some non-commutative
$4$-dimensional 
space.  A very nice, and not widely known, explanation of the $2\times 2 $ matrix
form was given by 
Vokos, Zumino and Wess \cite{Vokos}, who showed that the internal relations are
preserved under 
matrix multiplication in the following sense:

\newtheorem{theorem}{Theorem}

\newtheorem{vokos}{}
%\begin{vokos} If $U$ is as above  the entries of $U^n$ satisfy the
%relations $(R)$ 
%with $q$ substituted by $q^n$ for all positive integers $n$.
%\label{vokos1}
%\end{vokos}

\vspace{0.3cm}
\noindent 1.  {\em If $U$ is as above  the entries of $U^n$ satisfy the
relations $(R)$ 
with $q$ substituted by $q^n$ for all positive integers $n$.}
\vspace{0.3cm}

This result can be extended to all integers by making a minor modification.  
The quantum determinant $D_q=ad-qbc$ is central in $M_q(2)$. Formally adjoining 
a new generator $D_q^{-1}$, which commutes with $a,b,c,d$ and satisfies the 
relations $(ad-qbc)D_q^{-1}= 1$, gives rise to an algebra, 
denoted $GL_q(2)$, for which $U$ has the matrix inverse: 
\[ 
U^{-1}=D_q^{-1}\left( \begin{array}{cc}d &-q^{-1}b\\ -qc&a\end{array} \right). 
\] 

\vspace{0.3cm}
\noindent 2.  {\em If $U$ is as above with entries belonging to $GL_q(2)$, the
entries
of $U^n$ satisfy the 
relations $(R)$  with $q$ substituted by $q^n$
for all integers $n$.}
\vspace{0.3cm}

Finally suppose that  
\[ 
U'=\left( \begin{array}{cc}a' &b'\\ c'&d'\end{array} \right) 
\] 
is a second matrix with entries satisfying the same relations as those of $U$, 
denoted $(R')$, and commuting with the  entries of $U$. Make $U$ and $U'$ 
invertible by adjoining the generators $D_q^{-1}$ and $(D')_q^{-1}$, as above. 
Thus the entries of $U$ and $U'$ belong to the quotient of the  free algebra 
$k\left\{a,b,c,d, D_q^{-1}, a',b',c',d', (D')_q^{-1}\right\}$ by the ideal 
generated by $(R)$, $(R')$, commutativity relations between primed and unprimed 
generators, commutativity relations between $D_q^{-1}$, $(D')_q^{-1}$ and all 
other generators, and the relations $(ad-qbc)D_q^{-1}= 1$ and 
$(a'd'-qb'c')(D')_q^{-1}= 1$. 

\vspace{0.3cm}
\noindent 3.  {\em If $U$ and $U'$ are as above, the entries of $U^nU'^n$
satisfy the 
relations $(R)$ with $q$ replaced by $q^n$ for all integers $n$. }
\vspace{0.3cm}

These multiplicative properties of $2\times 2 $  quantum  matrices are very
striking, and provide the 
inspiration for the definition of quantum matrix pairs which follows.

Suppose, then, that we have two invertible matrices 
\[ 
U_1= \left( \begin{array}{cc}a &b\\ c&d\end{array} \right),\,\,\,\, U_2=\left( 
\begin{array}{cc} a' &b'\\ c'&d'\end{array} \right), 
\label{U1U2} 
\] 
 whose entries take values in a non-commutative algebra $\cal A$ over $k$, and
satisfy  
relations of two different types: internal relations $(I)$,  involving the 
entries of one matrix only and having the same structure for both matrices,  
and mutual relations $(M)$, involving the entries of both matrices at the same 
time. Both types 
of commutation relation may involve $q$, as well as possibly other scalar 
parameters. 
 
\newtheorem{defn}{Definition} 
\begin{defn} 
If the relations $(I)$ and $(M)$ for the matrices $U_1$ and  
$U_2$ are such that 
\begin{itemize} 
 \item[a)]  
\[ 
U_1U_2=qU_2U_1, 
%\label{fund} 
\] 
for some invertible $q\in k$, $q\neq 1$, where $q$ acts by scalar multiplication
on the right-hand-side, 
and  
\item[b)] 
the entries of $U_1^nU_2^m$ obey internal relations with the same structure  
$(I)$, for all integers $n$ and $m$,  up to possible substitution of 
scalar 
parameters,  
\end{itemize} 
we call $(U_1,U_2)$ a {\bf quantum matrix pair.}  
 
If condition $a)$ holds, but   
condition $b)$ only 
holds for $U_1^nU_2^n$, 
for all  integers $n$, we call $(U_1,U_2)$ a  {\bf restricted quantum 
matrix 
pair.} 
\label{qmpdef} 
\end{defn} 
 
Requirement $a)$ is a natural generalization of the quantum plane relation
\rref{qplane}, with $x$ and $y$ being replaced
by the $2\times 2$ matrices $U_1$ and $U_2$. Requirement $b)$, in either
its restricted or
unrestricted form, is analogous to the multiplicative  properties of quantum
matrices described above.

We will proceed to give examples of quantum matrix pairs, both restricted and
unrestricted.
These examples, based on our previous work \cite{2+1}, are relatively simple in
that both matrices are 
upper triangular, but still display novel features not found in diagonal
examples. In particular, the 
internal relations are not standard Weyl-type relations like \rref{qplane},
but involve three matrix 
entries simultaneously. We expect that quantum matrix pairs  for groups other
than $SL(2,R)$ can be 
found and will have an important role to play in the context of Chern-Simons
theory. 

In any case, quantum matrix pairs constitute an interesting
algebraic structure, worthy of study in its own right. The fact that they closely
resemble quantum 
groups could lead to novel insights and perspectives on the latter subject. Here
we should point out that 
there is a similar construction which arises in the context of Majid's {\em
braided matrices}  
\cite[Section 10.3]{maj}. He gives examples of pairs of $2\times 2$ matrices with
internal and mutual 
relations between their entries, such that the product of the first matrix with
the second obeys the same 
internal relations (but not the other way round). By comparison the internal
relations in our examples 
seem very robust, since they carry over to the product of the original matrices
in either order, as well as 
to other monomials in the two matrices or their inverses. Furthermore, two
different products of the 
original matrices may themselves form a new quantum matrix pair.

Perhaps the most intriguing feature about our examples is that their geometric
origin reemerges from 
the algebra, when we find an action of the modular group on spaces of quantum
matrix pairs.  A fuller 
understanding of the role of quantum matrix pairs in physical models will
undoubtedly involve notions  
of noncommutative {\em non-local} geometry. 

Our material is organized as follows. In section 2 the three types of example of 
quantum matrix pairs 
are given, and their individual features are discusssed. In section 3  several
ways are described  of
obtaining new quantum matrix pairs from a given one for the examples of section
2, which then leads 
to an action of the modular group on spaces of
quantum matrix pairs. Section 4 contains some comments.

\section{Examples of quantum matrix pairs \label{sec2}}

As pointed out in the introduction, the holonomy matrices are not gauge
invariant. Under gauge 
transformations they transform by simultaneous conjugation with an element of
$SL(2,R)$. 
Classically this means that the matrices can be simultaneously diagonalized, when
they are diagonalizable,
but there is also a sector where both matrices are (upper) triangular in form. We
will concentrate on the 
corresponding sector in the quantum theory, since it has much more interesting
behaviour than the 
diagonal case.

Thus from now on we will take both matrices $U_1$ and $U_2$ to be of
upper triangular form and will denote their entries as follows:
\be
U_i=\left( \begin{array}{cc} \a_i &\b_i\\ 0&\g_i\end{array}
\right), \,\,\,\,i=1,2
\label{Ui}
\ee
Also, for brevity, we will henceforth adopt the following convention:  when
the index $i$ appears in a statement, $i=1,2$ is understood. Since we
require both matrices $U_i$ to be invertible, we take $\a_i$ and $\g_i$ to
be invertible, which is formally achieved by adjoining new generators
$\a_i^{-1}$, $\g_i^{-1}$, and corresponding relations, to the free algebra
generated by $\a_i$, $\b_i$ and $\g_i$. The inverse of $U_i$ is then given by
\be
U_i^{-1}= \left( \begin{array}{cc} \a_i^{-1} & -\a_i^{-1}\b_i \g_i^{-1} \\
0 &  \g_i^{-1}
\end{array}\right).
\label{U-1}
\ee

In accordance with the definition of a quantum matrix pair, we must specify
internal (I) and
mutual (M) relations between the generators. These we will subdivide
further into diagonal (D) relations when 
they involve only diagonal entries, and
non-diagonal (ND) relations, when they also involve non-diagonal entries.
We will list various options for the relations below, and different
combinations of these options will furnish the three types of example we
want to present.

For the internal diagonal relations there are two choices, namely
\bea
{\rm (ID1)}\quad \a_i\g_i&=&\g_i\a_i=1\quad  \label{ID1}\\
{\rm (ID2)}\quad \a_i\g_i&=&\g_i\a_i\quad  \nonumber
\label{ID2}
\eea
The first choice implies $\g_i=\a_i^{-1}$, whereas the second 
merely requires the diagonal entries to commute. For the
non-diagonal internal
relations there are also two choices: 
\bea
{\rm (IND1)}\quad \a_i\b_i&=&\b_i\g_i,\quad 
\label{IND1}\\
{\rm (IND2)}\quad  \a_i\b_i&=&r\b_i\g_i,\quad  
\label{IND2}
\eea
where in the second choice 
$r$ is an invertible element of $k$, $r\neq 1$. Of course, $r$ and $q$ need not
be independent. For instance, they may be equal, or $r$ may be a power of $q$.

The mutual diagonal relations appear in the following two forms: 
\[
{\rm (MD1)}\quad  \a_1\a_2=q\a_2\a_1,   \label{MD1}
\]
\[
{\rm (MD2)}\quad \left\{ \begin{array}{lcr}\a_1\a_2&=&q\a_2\a_1,  \\
 \a_1\g_2&=&q^{-1}\g_2\a_1,   \\
 \a_2\g_1&=&q\g_1\a_2,   \\
 \g_1\g_2&=&q\g_2\g_1. \end{array} \right. 
\label{MD2}
\]
where $q\in k$ is the parameter for the fundamental $q$-commutation relation
\rref{fund}. In practice these two choices are the same, since we will
always combine (MD1) with (ID1), which implies (MD2), according to
Proposition \ref{ID1MD1->MD2} below.

Finally, we restrict ourselves to a single choice of mutual non-diagonal
relations:
\be
{\rm (MND)} \quad \left\{ \begin{array}{lcr}\a_1\b_2&=&q\b_2\g_1,  \\
 \b_1\g_2&=&q\a_2\b_1.  \end{array} \right.
\label{MND}
\ee

Before giving the examples, we need to derive a few results, starting with a
simple but important proposition,
which underlies all the subsequent calculations.

\newtheorem{prop}{Proposition}
\begin{prop} Let $\cal A$ be an algebra over the field $k$, and let 
$\a$, $\b$ and $\g$ be elements of $\cal A$, with $\a$ and $\g$
invertible. Let $q$ and $r$ be invertible elements of $k$. Then
\begin{itemize}
\item[a)] $\a\g=q\g\a\Rightarrow \a^n\g^m=q^{mn}\g^m\a^n,\,\,\forall 
n,m\in {\bf Z}$
\item[b)] $\a\b=r\b\g
\Rightarrow \a^n\b=r^n\b\g^n,\,\, \forall n\in {\bf Z}$.
\end{itemize}
\label{basic}
\end{prop}

\noindent {\em Proof:} a) This is trivial for $n,m\geq 0$, and follows 
from the relations $\a^{-1} \g=q^{-1}\g\a^{-1}$,
$\a\g^{-1}=q^{-1}\g^{-1}\a$ and
$\a^{-1}\g^{-1}=q\g^{-1}\a^{-1}$ when $n$ or $m$ or both are negative. 
b) This is trivial for  $n\geq 0$, and follows from 
$\a^{-1}\b=r^{-1}\b\g^{-1}$ when $n$ is negative. \hfill $\Box$
\vspace{1ex}

\noindent An immediate corollary of part a) is the result mentioned
after \rref{MD2}. 

\begin{prop}
\[
{\rm (ID1)}\wedge{\rm  (MD1)}\Longrightarrow {\rm (MD2)}
\]
\label{ID1MD1->MD2}
\end{prop}

The first requirement for $(U_1,U_2)$ to constitute a quantum matrix pair is the
relation a)
of the definition. This will be guaranteed in all examples by the mutual
relations (MD2) and (MND):
\be
U_1U_2=\left( \begin{array}{cc}\a_1\a_2 &\a_1\b_2+\b_1\g_2\\ 
0&\g_1\g_2\end{array}
\right)
=q\left( \begin{array}{cc}\a_2\a_1 &\a_2\b_1+\b_2\g_1\\ 
0&\g_2\g_1\end{array}
\right)=
qU_2U_1.
\label{fundpf}
\ee
To analyse the requirement b) of the definition we need expressions for
powers of the matrices $U_i$. These depend on the choice of internal
non-diagonal relations.
\begin{prop}
\beas
{\rm a)} & {\rm (IND1)}\Longrightarrow U_i^n\stackrel{\rm def}{=}
\left( \begin{array}{cc}\a_i(n)&\b_i(n)\\ 
0&\g_i(n)\end{array}\right)& =
\left(\begin{array}{cc}\a_i^n&n\b_i\g_i^{n-1}\\ 
0&\g_i^n\end{array}\right),\quad \forall n\in {\bf Z} 
\label{Uin1}\\
{\rm b)} & {\rm (IND2)}\Longrightarrow U_i^n\stackrel{\rm def}{=}
\left( \begin{array}{cc}\a_i(n)&\b_i(n)\\ 
0&\g_i(n)\end{array}\right)& =
\left(\begin{array}{cc}\a_i^n&\bar{n}_r \b_i\g_i^{n-1}\\ 
0&\g_i^n\end{array}\right),\quad \forall n\in {\bf Z}
\label{Uin2}
\eeas
where $\bar{n}_r=(1-r^n)/(1-r)$.
\label{Uin}
\end{prop}
\noindent {\em Proof:} For $n\geq 0$ the formulae are proved by 
induction. The induction step uses the equalities:
\bea
\a_i^n\b_i+n\b_i\g_i^n&=&(n+1)\b_i\g_i^n \nonumber\\
\a_i^n\b_i+\bar{n}_r \b_i\g_i^n&=&(r^n+\bar{n}_r ) \b_i\g_i^n =
(\overline{n+1})_r\b_i\g_i^n \nonumber
\eea
for (IND1) and (IND2) respectively. These equalities follow from 
Proposition
\ref{basic}. 

For negative $n$, set $n=-p$ with $p$ positive. $U_i^{-1}$ has the 
internal relations:
\bea
\a_i(-1)\b_i(-1)&=& \a_i^{-1} (-\a_i^{-1}  \b_i \g_i^{-1})=
(-\a_i^{-1}\b_i \g_i^{-1}) \g_i^{-1}=\b_i(-1)\g_i(-1),\nonumber\\
\a_i(-1)\b_i(-1)&=& \a_i^{-1} (-\a_i^{-1}  \b_i \g_i^{-1})=
r^{-1}(-\a_i^{-1}\b_i \g_i^{-1})
\g_i^{-1}=r^{-1}\b_i(-1)\g_i(-1)\nonumber
\eea
for (IND1) and (IND2) respectively, using \rref{U-1} and Proposition
\ref{basic}. Since $U^n=(U^{-1})^p$ one derives:
\bea
\b_i(n)&=&p\b_i(-1)\g_i(-1)^{p-1}\nonumber\\
&=& p (-\a_i^{-1}  \b_i\g_i^{-1})\g_i(-1)^{-(p-1)}\nonumber\\
&=&-p\b_i\g_i^{-p-1}=n\b_i\g_i^{n-1}\nonumber
\eea
and
\bea
\b_i(n)&=&\bar{p}_{r^{-1}}\b_i(-1)\g_i(-1)^{p-1}\nonumber\\
&=& \bar{p}_{r^{-1}}(-\a_i^{-1}\b_i\g_i^{-1})\g_i(-1)^{-(p-1)}\nonumber\\
&=&-r^{-1}\bar{p}_{r^{-1}}\b_i\g_i^{-p-1}
=\bar{n}_r\b_i\g_i^{n-1}\nonumber
\eea
for (IND1) and (IND2) respectively. The final equality is the 
calculation:
\[
-r^{-1} \bar{p}_{r^{-1}}=- r^{-1}\frac{(1- (r^{-1})^p)}{1- r^{-1}}=
-\frac{(1-r^{-p})}{r-1}=\frac{1-r^n}{1-r}=\bar{n}_r. \quad \Box 
\]
As a corollary we obtain formulae for the entries of $U_1^nU_2^m$ for
$n,m\in {\bf Z}$. We remark that, in view of the $q$-commutation relation
\rref{fund} and Proposition \ref{basic}, any word in $U_1$, $U_2$ and their
inverses is proportional to $U_1^nU_2^m$ for some $n,m\in {\bf Z}$. Setting
\be
U_1^nU_2^m=\left( \begin{array}{cc}\a(n,m)&\b(n,m)\\ 
0&\g(n,m)\end{array}\right),
\label{Uinm,code}
\ee
we obtain the following formulae:
\bea
& \a(n,m)=&\a_1^n\a_2^m\nonumber\\
  & \g(n,m)=&\g_1^n\g_2^m\nonumber\\
{\rm (IND1)} \Longrightarrow & \b(n,m)=&
m\a_1^n\b_2\g_2^{m-1} + n\b_1\g_1^{n-1}\g_2^m\nonumber\\
{\rm (IND2)} \Longrightarrow & \b(n,m)=&
\bar{m}_r\a_1^n\b_2\g_2^{m-1} + \bar{n}_r\b_1\g_1^{n-1}\g_2^m\nonumber
\label{Uinm,entries}
\eea

With these preliminary calculations out of the way, we are in a position to
present our three types of example, and prove that they are quantum matrix pairs.
We do this
in the form of a theorem.
\begin{theorem}
Matrix pairs of the form  {\em \rref{Ui}}, satisfying internal and mutual
relations as set out in the table below, 
give rise to three types of quantum matrix pairs, of which the third is a
restricted quantum matrix pair.
{\em
\begin{center}
\begin{tabular}{|l|l|l|} \hline
Example & Internal relations & Mutual relations \\ \hline
Type I &  (ID1), (IND1) & (MD1), (MND) \\ \hline
Type II &  (ID2), (IND1) & (MD2), (MND) \\ \hline
Type III &  (ID2), (IND2) & (MD2), (MND) \\ \hline
\end{tabular}
\end{center}
}
For type I, $U_1^nU_2^m$, for all $n,m\in {\bf Z}$, has internal relations:
\bea
\a(n,m)\g(n,m)&=& \g(n,m)\a(n,m) =\,\, q^{nm}
\label{IA}\\ 
\a(n,m)\b(n,m)&=& \b(n,m)\g(n,m).  
\label{IB}
\eea
For type II, $U_1^nU_2^m$, for all $n,m\in {\bf Z}$, has internal relations:
\bea
\a(n,m)\g(n,m)&=& \g(n,m)\a(n,m) 
\label{IIA}\\ 
\a(n,m)\b(n,m)&=& \b(n,m)\g(n,m).  
\label{IIB}
\eea
For type III,  $U_1^nU_2^n$, for all
$n\in {\bf Z}$, has internal relations:
\bea
\a(n,n)\g(n,n)&=& \g(n,n)\a(n,n)
\label{IIIA}\\ 
\a(n,n)\b(n,n)&=& r^n\b(n,n)\g(n,n).  
\label{IIIB}
\eea
\label{expsthm}
\end{theorem}

In each case these internal relations have the same structure as those of
the corresponding $U_i$, with $1$ in \rref{ID1} replaced with $ q^{nm}$ in
\rref {IA} for type I, and $r$ in \rref{IND2} replaced with $r^n$ in \rref
{IIIB} for type III.
Whilst these three types of example are very similar, each of them exhibits
some 
special feature 
distinguishing it from the others. In the type I case, both $U_1$ and $U_2$
have determinant 1, but mixed products of the $U_i$ have non-unit determinant,
as a result of the non-commutativity of the algebra. In the type II case,
despite the noncommutativity, the property of having commuting diagonal
entries propagates to all products of the $U_i$. Finally, the type III case has
internal relations involving a parameter, a feature which is reminiscent of
quantum groups. We remark that this parameter cannot simply be removed by a
rescaling $\a_i\rightarrow r^{1/2}\a_i$, $\g_i\rightarrow r^{-1/2}\g_i$,
since the (MND) relations are not preserved  under these replacements. The
parameter propagates to powers of $U_i$ in a manner again reminiscent of
quantum groups (cf. the second Vokos et al result in the introduction).

\noindent {\em Proof of Theorem \ref{expsthm}:} Because of
Proposition \ref{ID1MD1->MD2} and equation \rref{fundpf}, all three
types satisfy requirement a) of the definition of a quantum matrix
pair.  All that remains is to show that the relations
\rref{IA}-\rref{IIIB} hold.  Equations \rref{IA}, \rref{IIA} and
\rref{IIIA} follow from 
\begin{eqnarray*}
(\a_1^n\a_2^m)(\g_1^n\g_2^m)&=&q^{mn}\a_1^n\g_1^n\a_2^m\g_2^m (=
q^{mn})\\ &=& q^{mn}\g_1^n\a_1^n\g_2^m\a_2^m=
(\g_1^n\g_2^m)(\a_1^n\a_2^m), \label{Apf} 
\end{eqnarray*} 
where the
equality between brackets holds for the type I case, and using
Proposition \ref{basic} a) in the first and last equalities.
\rref{IB} and \rref{IIB} follow from the calculation
\begin{eqnarray*}
\lefteqn{(\a_1^n\a_2^m)(m\a_1^n\b_2\g_2^{m-1}+n\b_1\g_1^{n-1}\g_2^m)}
\nonumber\\ 
&=& mq^{-mn}\a_1^n(\a_1^n\a_2^m)\b_2\g_2^{m-1} +
nq^{-m}\b_1(\g_1^n\g_2^m)\g_1^{n-1}\g_2^m\nonumber\\ 
&=&mq^{-mn}q^n\a_1^n\b_2(\g_1^n\g_2^m)\g_2^{m-1} +
nq^{-m}q^{-m(n-1)}\b_1\g_1^{n-1}(\g_1^n\g_2^m)\g_2^m\nonumber\\ 
&=&mq^{-mn}q^nq^{n(m-1)}\a_1^n\b_2\g_2^{m-1}(\g_1^n\g_2^m) +
nq^{-m}q^{-m(n-1)}q^{nm}\b_1\g_1^{n-1}\g_2^m(\g_1^n\g_2^m)\nonumber\\
&=& (m\a_1^n\b_2\g_2^{m-1}+n\b_1\g_1^{n-1}\g_2^m)(\g_1^n\g_2^m),
\label{IBIIBpf} 
\end{eqnarray*} 
where we have repeatedly used
Proposition \ref{basic}, as well as Proposition \ref{ID1MD1->MD2}
for the type I case. For \rref{IIIB}, the calculation is slightly
modified due to the presence of the parameter $r$ in the (IND2) 
relation:  
\beas
\lefteqn{(\a_1^n\a_2^n)(\a_1^n\b_2\g_2^{n-1}+\b_1\g_1^{n-1}\g_2^n)}
\nonumber\\ 
&=& q^{-n^2}\a_1^n(\a_1^n\a_2^n)\b_2\g_2^{n-1} +
r^nq^{-n}\b_1(\g_1^n\g_2^n)\g_1^{n-1}\g_2^n\nonumber\\ 
&=&q^{-n(n-1)}r^n\a_1^n\b_2(\g_1^n\g_2^n)\g_2^{n-1} +
r^nq^{-n^2}\b_1\g_1^{n-1}(\g_1^n\g_2^n)\g_2^n \nonumber\\ 
&=&r^n(\a_1^n\b_2\g_2^{n-1}+\b_1\g_1^{n-1}\g_2^n)(\g_1^n\g_2^n) 
\label{IIIBpf} 
\eeas 
where we have omitted the factor $n$ in
$\b(n,n)$. 
\hfill $\Box$

\section{Generating new quantum matrix pairs \label{sec3}}

The examples of the previous section showed 
how the internal relations are preserved under multiplication of the
matrices belonging to a quantum matrix pairs.
However there is another aspect to these examples. By taking two different
products of
the $U_i$, in some circumstances it is possible to generate a new quantum matrix
pair of
the same or a similar type, as we will see in this section. Furthermore the
transformations amongst  quantum matrix pairs may preserve the type, so that we
can
regard them as acting on the space of all quantum matrix pairs of a certain type.
We show
how this can give rise to representations of a discrete group, namely
$SL(2,{\bf Z})$ (the modular group), on the space of quantum matrix pairs of 
type I and
type II.

We start with a trivial first result in this direction.
\begin{prop} Let $(U_1,U_2)$ be a quantum matrix pair of any of the three types
described in
the previous section. Then the pair $(\tilde{U} _1, \tilde{U} _2)$ with
entries $\tilde{\a}_i=\a_i$, $\tilde{\g}_i=\g_i$ and $\tilde{\b}_i=c_i\b_i$,
where $c_i\in k$ are arbitrary constants, is a new quantum matrix pair of the
same type as
the original pair.
\label{rescalingprop}
\end{prop}
\noindent {\em Proof:} The non-diagonal relations \rref{IND1}, \rref{IND2}
and \rref{MND} are linear in $\b_1$ or $\b_2$.
\hfill $\Box$
\vspace{1ex}

\noindent
The main result of this section is stated in the following theorem:
\begin{theorem} 
\begin{itemize}
\item[a)] Let $(U_1,U_2)$ be a quantum matrix pair of type I. Then 

\noindent $(q^{-nm/2}U_1^nU_2^m, q^{-st/2}U_1^sU_2^t)$ is a quantum matrix
pair
of type I, with $q$
replaced with $q^{nt-ms}$, for all $n,m,s,t\in {\bf Z}$.
\item[b)] Let $(U_1,U_2)$ be a quantum matrix pair of type II. Then 
$( U_1^nU_2^m, U_1^sU_2^t)$ is a quantum matrix pair of type II, with $q$
replaced with
$q^{nt-ms}$, for all $n,m,s,t\in {\bf Z}$.
\item[c)] Let $(U_1,U_2)$ be a quantum matrix pair of type III. Then 
$( U_1^n,U_2^n)$ is a quantum matrix pair of type III, with $q$ replaced with 
$q^{n^2}$ and $r$ replaced with $r^n$, for all $n \in {\bf Z}$.
\end{itemize}
\label{qmp->qmp}
\end{theorem}
\noindent {\em Proof:} (The statements in the proof hold for all $n,m,s,t\in {\bf
Z}$.) First we prove the internal relations. For a), the
(ID1) relations follow from \rref{IA}, since this equation implies that
$q^{-nm/2}\a(n,m)$ and $q^{-nm/2}\g(n,m)$ are each other's inverses.
The (IND1) relations follow from \rref{IB}, since
all matrix entries are multiplied by the same factor. For b), the internal
relations are equations \rref{IIA}, \rref{IIB} of the previous section. For
c), using the notation and result of Proposition \ref{Uin} b),
$a_i\g_i=\g_i\a_i$ (ID2) for $U_i$ implies
$\a_i(n)\g_i(n)=\g_i(n)\a_i(n)$, and
$\a_i\b_i=r\b_i\g_i$ (IND2) for $U_i$, together with Propostion \ref{basic},
implies $\a_i(n)\b_i(n)=r^n\b_i(n)\g_i(n) )$.

To simplify the proof of the mutual relations , we use the relation
\[
( U_1^nU_2^m)( U_1^sU_2^t)= q^{nt-ms}( U_1^sU_2^t)
( U_1^nU_2^m),\nonumber
\]
which follows from the $q$-commutation relation \rref{fund}  satisfied by
$U_i$ and Proposition \ref{basic}. This implies the equations
\bea
\a(n,m)\a(s,t)&=&q^{nt-ms}\a(s,t)\a(n,m) \label{new11}\\
\g(n,m)\g(s,t)&=&q^{nt-ms}\g(s,t)\g(n,m) \label{new22}\\
\a(n,m)\b(s,t) + \b(n,m)\g(s,t)&=&q^{nt-ms}(\a(s,t)\b(n,m) +
\b(s,t)\g(n,m)
\label{new12}
\eea
where we are using the notation of
\rref{Uinm,code}. 

Now, starting with b), the first and fourth (MD2)
relations, with $q$ replaced with $q^{nt-ms}$, are equations \rref{new11}
and \rref{new22}, and the second and third (MD2) relations
\bea
\a(n,m)\g(s,t)&=&q^{-(nt-ms)}\g(s,t)\a(n,m) \nonumber\\
\a(s,t)\g(n.m)&=&q^{nt-ms}\g(n,m)\a(s,t) \nonumber
\eea
follow from the second and third (MD2)
relations for $U_i$ and Proposition \ref{basic}. In view of \rref{new12}, it
is enough to show the first of the (MND) relations:
\bea
\a(n,m)\b(s,t)&=&(\a_1^n\a_2^m)(t\a_1^s\b_2\g_2^{t-1}
+s\b_1\g_1^{s-1}\g_2^t) \nonumber\\
&=& q^{-ms}t\a_1^s(\a_1^n\a_2^m)\b_2\g_2^{t-1}+
sq^{-m}\b_1(\g_1^n\g_2^m)\g_1^{s-1}\g_2^t \nonumber\\
&=& q^{-ms+n} t \a_1^s\b_2(\g_1^n\g_2^m)\g_2^{t-1} + q^{-ms}s
\b_1\g_1^{s-1}(\g_1^n\g_2^m)\g_2^t\nonumber\\
&=& q^{nt-ms}(t\a_1^s\b_2\g_2^{t-1}+s\b_1\g_1^{s-1}\g_2^t) )(\
g_1^n\g_2^m)\nonumber\\
&=&q^{nt-ms}\b(s,t)\g(n,m)\nonumber
\eea

For a), the (MD1) relation follows from \rref{new11}, and the (MND) relation
is proved as for case b), after multiplying the entries by the factors
$q^{-nm/2}$ or $q^{-st/2}$. 
For c), the (MD1) relations, with parameter $q^{n^2}$, are shown as for b),
after setting $t=n,\,m=s=0$. Again, in view of \rref{new12}, it is enough to
show the first of the (MND) relations:
\bea
\a_1(n)\b_2(n)&=&\a_1^n\bar{n}_r\b_2\g_2^{n-1}=
\bar{n}_rq^n\b_2\g_1^n\g_2^{n-1}\nonumber\\
&=& \bar{n}_rq^{n^2}\b_2\g_2^{n-1}\g_1^n= q^{n^2}\b_2(n)\g_1(n). \quad \Box
\nonumber
\eea

The modular group $SL(2,{\bf Z})$ has a presentation in terms of two
generators $S$ and $T$, with relations $S^4=(ST)^3=1$.\footnote{For a study of
the
modular group in the context of $(2+1)$ quantum gravity by Carlip and one of the
authors, 
see  \cite{qgmod}.} We
have a representation of  $SL(2,{\bf Z})$, if we can find automorphisms, $S$
and $T$, of a space $X$, which satisfy these relations.
There are
natural $SL(2,{\bf Z})$ representations associated with spaces of quantum matrix
pairs, as
the following theorem shows.

\begin{theorem}
Let ${\rm QMP}1$ and ${\rm QMP}2$ be the spaces of all quantum matrix pairs of
type I and II
respectively, with entries in an algebra $\cal A$, and with a fixed parameter
$q$. Then the following
definitions give rise to representations of $SL(2,{\bf Z})$ on ${\rm QMP}1$  and
${\rm QMP}2:$
\bea
{\rm QMP}1:& S(U_1,U_2)=(U_2,U_1^{-1}), & T(U_1,U_2)=(
q^{-1/2}U_1U_2,U_2)\nonumber\\
{\rm QMP}2:& S(U_1,U_2)=(U_2,U_1^{-1}), &
T(U_1,U_2)=( U_1U_2,U_2)\nonumber
\eea
\label{modular}
\end{theorem}
\noindent {\em Proof:} From a) and b) of the previous theorem, the
transformations $S$ and $T$ map type I or II quantum matrix pairs into quantum
matrix pairs of the same type, and with the same $q$ parameter. 
$S^2(U_1,U_2)=S(U_2,U_1^{-1})= (U_1^{-1}, U_2^{-1})$, and thus
$S^4(U_1,U_2)=(U_1,U_2)$. For the type I case, the second relation
$(ST)^3=1$ is proved as follows. First:
\[
(ST)(U_1,U_2)=S(q^{-1/2}U_1U_2,U_2)=(U_2, q^{1/2}U_2^{-1}
U_1^{-1}).
\]
Thus
\bea
(ST)^3(U_1,U_2)&=&(ST)^2(U_2, q^{1/2}U_2^{-1} U_1^{-1}) \nonumber\\
&=& (ST)( q^{1/2}U_2^{-1} U_1^{-1}, q^{1/2}(q^{-1/2}U_1U_2)U_2^{-1})
\nonumber \\
&=& (ST)( q^{1/2}U_2^{-1} U_1^{-1}, U_1)\nonumber\\
&=&(U_1,q^{1/2}U_1^{-1}(q^{-1/2}U_1U_2))= (U_1,U_2). 
\nonumber
\eea
For the type II case, set $q=1$ in this calculation. \hfill $\Box$

\section{Final comments \label{sec4}}

Quantum matrix pairs combine the preservation of internal relations under
multiplication, a quantum-group-like feature, with the fundamental
$q$-commutation relation which holds between the two matrices.  We have presented
three types of example of this construction, all involving upper-triangular
matrices, but with slightly differing features.

It is interesting to make some comparisons between quantum matrix pairs and
quantum groups. The internal relations in our examples of quantum matrix pairs
differ in structure from the Weyl-type q-commutation relations normally found in
quantum groups, as they involve three matrix elements at the same time.  Related
to this is the fact that the entries of each matrix do not commute in the limit
$q\rightarrow 1$, which also distinguishes them from Majid's braided matrices
\cite{maj}.  Nonetheless, when the internal relations depend on a parameter, as
in the third type of example, quantum integers with that parameter appear in the
powers of the matrices, which is a feature strongly reminiscent of quantum
groups.

An obvious question  for further study is to see whether other examples can be
found, e.g. $2 \times 2$ matrices but not of triangular form, or examples
involving other groups.

According to theorem \ref{qmp->qmp}, not only do products of powers of the 
matrices have the
same structure of internal relations, but taking two different products gives
rise to new quantum matrix pairs of the same type.  This shows that, in a sense,
it is the whole quantum matrix pair structure, rather than just the internal
relations, which is preserved under multiplication in these examples.

It is striking that the action of the modular group on pairs of commuting
matrices extends to quantum matrix pairs. This reveals that the construction,
which could be taken on a  purely algebraic level, actually has a geometric
interpretation as well. In future work we hope to arrive at a deeper
understanding of quantum matrix pairs in terms of non-local non-commutative
geometry.

\section*{Acknowledgments}
This work was supported in part by the European Commission HCM programme
CHRX-CT93-0362, the European Commission TMR programme
ERBFMRX-CT96-0045, INFN Iniziativa Specifica FI41, and the
Programa de Financiamento Plurianual of the Funda\c{c}\~{a}o 
para a Ci\^{e}ncia e a Tecnologia (FCT).

\end{document}